\documentclass[12pt,a4paper, twoside]{amsart}

\usepackage[all]{xypic}
\usepackage{latexsym, amsmath, amscd, amssymb, amsthm, epsfig, times}
\usepackage{graphicx}

\usepackage{color}
\definecolor{maus}{rgb}{0.2,0.2,0.2}

\newtheorem{thm}{Theorem}

\newtheorem{lem}[thm]{Lemma}
\newtheorem{cor}[thm]{Corollary}

\theoremstyle{definition}
\newtheorem{defn}[thm]{Definition}

\newtheorem*{ex*}{Example}

\newcommand{\thmref}[1]{Theorem~\ref{#1}}

\newcommand{\R}{\mathbb{ R}}

\newcommand{\tb}{\operatorname{tb}}
\newcommand{\rot}{\operatorname{rot}}
\newcommand{\sel}{\operatorname{sl}}

\newcommand{\p}{\partial}

\newcommand{\eing}[1]{\big|_{#1}}

\title{Legendrian ribbons in overtwisted contact structures}

\author{S. Baader, K. Cieliebak, T. Vogel}

\thanks{KC is partially supported by DFG grant CI 45/2-2}

\address{S. Baader, Dept. Mathematik der ETH Z{\"u}rich,
  R{\"a}mistr. 101, 8092 Z{\"u}rich, Switzerland} 
\email{sebastian.baader@math.ethz.ch}

\address{K. Cieliebak, Mathematisches~Institut der LMU
M\"unchen, Theresienstr.~39, 80333 M\"unchen, Germany}
\email{kai@math.lmu.de}

\address{T. Vogel, Mathematisches~Institut der LMU
M\"unchen, Theresienstr.~39, 80333 M\"unchen, Germany}
\email{tvogel@math.lmu.de}

\date{\today}

\begin{document}
\color{maus}

\begin{abstract}
We show that a null--homologous transverse knot $K$ in the complement
of an overtwisted disk in a contact 3--manifold is the boundary of a
Legendrian ribbon if and only if it possesses a Seifert surface $S$
such that the self--linking number of $K$ with respect to $S$ satisfies
$\sel(K,S)=-\chi(S)$. In particular, every null--homologous topological
knot type in an overtwisted contact manifold can be represented by the 
boundary of a Legendrian ribbon. Finally, we show that a contact
structure is tight if and only if every Legendrian ribbon minimizes
genus in its relative homology class. 
\end{abstract}

\maketitle

\section{Introduction}

In this note, $(M,\xi)$ will always denote a cooriented contact 3--manifold, i.e.~$\xi=\ker\lambda$ for a 1--form $\lambda$ on $M$ satisfying $\lambda\wedge d\lambda\neq 0$. An oriented knot $K \subset M$ is called {\em(positively) transverse} if $\lambda(\dot\gamma)>0$ for a positive parametrization $\gamma:S^1\to K$. We will always assume that $K$ is null--homologous. Then $K$ possesses a {\em Seifert surface}, i.e.~an embedded oriented connected surface $\Sigma\subset M$ with boundary $\p\Sigma=K$. Given $K$ and $\Sigma$, we choose a nowhere vanishing section $v$ of $\xi\eing{\Sigma}$ and use $v$ to push $K$ away from itself. The resulting knot is denoted by $K'$. 

\begin{defn}
The {\em self--linking number} $\sel(K,\Sigma)$ is defined as the algebraic intersection number of $K'$ and $\Sigma$. 
\end{defn}

The self--linking number $\sel(K,\Sigma)$ depends only on $K$ and $[\Sigma]\in H_2(M,K)$ and is independent of the choice of $v$. Moreover it is always odd (since it has the same parity as the Euler characteristic of $\Sigma$). When it is clear which Seifert surface we use we will simply write $\sel(K)$. Obviously, isotopic transverse knots have the same self--linking number (with respect to Seifert surfaces carried along with the isotopy). 

A contact structure $\xi$ on a manifold $M$ is called \emph{overtwisted} if $M$ contains an overtwisted disc, i.e. an embedded disc whose boundary is a Legendrian unknot with Thurston--Bennequin number zero (see Figure~1). Otherwise the contact structure $\xi$ is called \emph{tight}. 

\begin{figure}[ht]
\center
\includegraphics{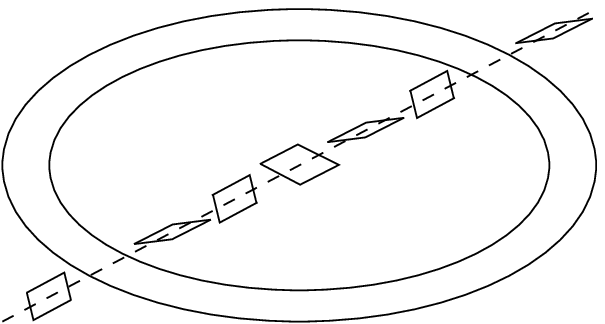}
\caption{}
\end{figure}

The following result was first proved for the standard contact structure on $\R^3$ by D.~Bennequin in \cite{Be} and then generalized to tight contact structures by Y.~Eliashberg in \cite{El}.  

\begin{thm}\label{t:bennequin}
If $(M,\xi)$ is a tight contact 3--manifold and $K$ is a transverse knot with Seifert surface $\Sigma$, then  
\begin{equation} \label{e:tb ungl}
\sel(K,\Sigma)\le-\chi(\Sigma).
\end{equation}
\end{thm}

\begin{defn}
A {\em Legendrian graph} in a contact 3--manifold $(M,\xi)$ is a trivalent embedded graph $\Gamma\subset M$ such that all edges are tangent to $\xi$. 
\end{defn}

To a Legendrian graph one can associate a transverse link type as follows. Choose a surface $\Sigma_\Gamma$ containing $\Gamma$ with smooth boundary such that $\Sigma_\Gamma$ is tangent to $\xi$ at every point of $\Gamma$.

We use the orientation of the contact structure to orient $\Sigma_\Gamma$. Then for $\Sigma_\Gamma$ sufficiently small (i.e.~after replacing $\Sigma_\Gamma$ by a sufficiently small neighbourhood of $\Gamma$ in $\Sigma_\Gamma$), its boundary $\partial \Sigma_\Gamma$ is a link all of whose components are positive transverse knots. The isotopy class of the resulting transverse link depends only on $\Gamma$. We refer to $\Sigma_\Gamma$ as a {\em Legendrian ribbon}.

\section{Legendrian Ribbons and overtwisted discs}

The following theorem is an analogue of a result explained in \cite{Dy} for overtwisted knots. Before stating it, recall that positive transverse knots are obtained from knots tangent to $\xi$ by choosing an oriented framing of $\xi$ such that one of the components of the framing is tangent to the knot and pushing the knot in the direction opposite to the second component of the framing. As described for example in in \cite{Et}, a positive transverse knot $\gamma$ can be stabilized to a knot $\hat{\gamma}$ with $\sel(\hat{\gamma})=\sel(\gamma)-2$.  

\begin{thm} \label{t:class}
Let $(M,\xi)$ be a contact 3--manifold and $D\subset M$ an overtwisted disc. If $L,K$ are two null--homologous transverse knots lying in the complement of $D$ such that $K,L$ represent the same topological knot type and $\sel(K)=\sel(L)$ (with respect to Seifert surfaces carried along with the isotopy), then $L$ and $K$ are isotopic as transverse knots.  
\end{thm}

\begin{proof}
Any two transverse knots $K,L$ representing the same topological knot types become transversely isotopic after sufficiently many stabilizations; this fact can be shown as the analogous statement for Legendrian knots, cf. \cite{Et}. (At this point the contact structures are allowed to be tight or overtwisted.) 

Of course, every stabilization changes the self--linking number of the knots and therefore this procedure by itself does not produce a transverse isotopy between the original knots $K,L$. However, when $K,L$ lie in the complement of a fixed overtwisted disc $D$, then one can neutralize each stabilization by pulling a segment of the knot over $D$. Indeed, let $\partial D$ be oriented in such a way that $\rot(\partial D)=1$, then its positive transversal push--off $\partial D^+$ has self--linking number one, according to the general formula
$$
\sel(\gamma^+)=\tb(\gamma)+\rot(\gamma),
$$
which holds for all null--homologous Legendrian knots $\gamma$ in a 3--manifold (see \cite{Be}). Here $\tb$ and $\rot$ denote the Thurston--Bennequin number and the rotation number, respectively. 

Given two transversal knots $\gamma, \gamma'$ in $M$ which lie in disjoint balls one can define the connected sum such that the self linking number of the resulting knot $\gamma\#\gamma'$ satisfies
$$
\sel(\gamma\#\gamma')=\sel(\gamma)+\sel(\gamma')+1.
$$
Let $D_{ot}$ be an overtwisted disc and $\gamma'$ the positive push off of the boundary $\partial D_{ot}$ with $\sel(\gamma')=+1$. If $\gamma$ is a positive transverse knot which is disjoint from $D_{ot}$, then $\sel(\gamma\#\gamma')=\sel(\gamma)+2$ and $\gamma\#\gamma'=\gamma$ as topological knot types. We call $\gamma\#\gamma'$ the destabilization of $\gamma$.  A stabilization of $\gamma\#\gamma'$ yields a knot which is isotopic to $\gamma$ as positive transverse knot. The isotopy is obtained by a push off of the isotopy constructed in Lemma~4.7 in \cite{Dy} (where a similar situation for Legendrian knots is considered). 

In order to construct a transverse isotopy between $K$ and $L$ we follow a procedure used in \cite{Dy} in the context of Legendrian knots.  We first isotope a segment of $K$ such that it coincides with a segment $\sigma$ of $L$. Then we stabilize both knots sufficiently often on the complements of $\sigma$ and we destabilize  $\sigma$ using $D$ in order to undo the stabilizations. 

If this is done sufficiently many times, then the complements of $\sigma$ become transversely isotopic while $K$ and $L$ still coincide along $\sigma$. Because we have never changed the transverse knot types this shows that $K,L$ are isotopic transverse knots.  
\end{proof}

The following lemma shows that transverse knots obtained from Legendrian ribbons realize equality in the Thurston--Bennequin inequality \eqref{e:tb ungl}. Y.~Kanda has given examples of topological knot types for which the equality in \eqref{e:tb ungl} cannot be realized by any transverse representative, cf. \cite{Ka}. In particular, such knots cannot be boundaries of Legendrian ribbons.  

\begin{lem} \label{l:euler}
Let $K=\partial\Sigma_\Gamma$ be the boundary of a Legendrian ribbon $\Sigma_\Gamma$. Then $\sel(K,\Sigma_\Gamma)=-\chi(\Sigma_\Gamma)$. 
\end{lem}

\begin{proof}
Let $w$ be a vector field along $K=\p\Sigma_\Gamma$ which is tangent to $\Sigma_\Gamma$ and points outwards. As $\Sigma_\Gamma$ is almost tangent to $\xi$, $w$ projects to a nonvanishing section $v$ of $\xi$ along $K$. The self--linking number of $K$ can be described as the obstruction to the existence of a nowhere vanishing extension over $\Sigma_\Gamma$ of $v$ as a section of $\xi$, or equivalently, of $w$ as a section of $T\Sigma_\Gamma$. But the latter obstruction equals minus the Euler characteristic, hence $\sel(K,\Sigma_\Gamma)=-\chi(\Sigma_\Gamma)$.  
\end{proof}

\begin{thm}\label{t:ribbon}
Let $K$ be a transverse knot in an overtwisted contact 3--manifold $(M,\xi)$ with Seifert surface $\Sigma$. Assume that $M\setminus K$ is still overtwisted and that $\sel(K,\Sigma)=-\chi(\Sigma)$. Then there is a Legendrian graph $\Gamma$ with $\partial\Sigma_\Gamma=K$. 
Moreover, $[\Sigma_\Gamma]=[\Sigma]\in H_2(M,K)$ and $\chi(\Sigma_\Gamma)=\chi(\Sigma)$.  
\end{thm}

\begin{proof}
We fix an overtwisted disc $D$ in the complement of $K$. Without loss of generality we may assume that $D$ and $\Sigma$ do not intersect. The following constructions can be carried out in the complement of $D$.  

Let $G\subset \Sigma$ be a trivalent graph such that $G$ is a deformation retract of $\Sigma$. After a $C^0$--small isotopy of $\Sigma$ we may assume that $\Sigma$ is tangent to $\xi$ at the vertices of $G$ such that the orientations of $\Sigma$ and of $\xi$ coincide at the vertices. Moreover, we may assume that the edges of $G$ are tangent to $\xi$ near the vertices.  

Next, we smoothly isotope each edge to a Legendrian curve contained in a small neighbourhood of the original edge, fixing it near the vertices, and pull along $\Sigma$ with the isotopy. After doing this, $\Sigma$ and $\xi$ need  not induce the same framing of the edges. However, the framings can be arranged to agree for each edge by either stabilizing the edge sufficiently often or by  sliding the edge sufficiently often over $D$, cf. \cite{Dy}. These operations correspond to taking a connected sum with a Legendrian unknot with Thurston--Bennequin number $-1$ or $+1$, respectively. The latter only exist in the presence of an overtwisted disc. An example of a Legendrian ribbon whose core curve is a Legendrian unknot with Thurston--Bennequin number $+1$ is shown in Figure~2. The two `horizontal' parts of the ribbon lie on the boundary of two parallel overtwisted discs there.

\begin{figure}[ht]
\center
\includegraphics{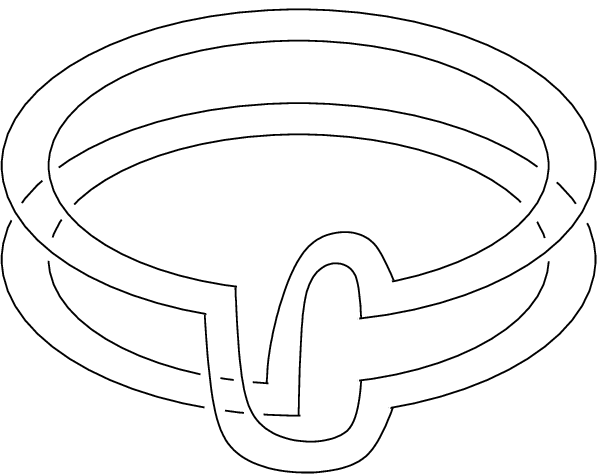}
\caption{}
\end{figure}

We denote the resulting Legendrian graph by $\Gamma$ and the resulting transverse knot by $L=\p\Sigma_\Gamma$.   

By construction, $K$ and $L$ are equivalent as topological knots, $\sel(K)=-\chi(\Sigma)=-\chi(\Sigma_\Gamma)=\sel(L)$ by Lemma~\ref{l:euler}, and the complement of $K\cup L$ is still overtwisted. So by \thmref{t:class} the knots $K$ and $L$ are transversely isotopic. This isotopy can be realized by an ambient contact isotopy (see~\cite{Et}), so after pulling along $\Sigma_\Gamma$ with this isotopy we may assume $\p\Sigma_\Gamma=K$. By construction we have $[\Sigma_\Gamma]=[\Sigma]\in H_2(M,K)$ and $\chi(\Sigma_\Gamma)=\chi(\Sigma)$, 
which finishes the proof of the theorem. 
\end{proof}

Denote by $g(S)$ the genus of a surface $S$. 

\begin{cor}
Let $(M,\xi)$ be an overtwisted contact manifold. Let $K\subset M$ be a null--homologous transverse knot in the complement of an overtwisted disc and $S$ be a Seifert surface for $K$ which minimizes the genus in its relative homology class.  

Then $K$ is the bondary of a Legendrian ribbon representing the class $[S]\in H_2(M,K)$ if and only if $\sel(K,S)\geq 2g(S)-1$. In particular, every null--homologous topological knot type in an overtwisted contact $3$--manifold can be represented by the boundary of a Legendrian ribbon.  
\end{cor}

\begin{proof}
Suppose first that $\sel(K,S)\geq 2g(S)-1$. By taking connected sums with null--homologous tori we can increase the genus of $S$ by any positive integer, without changing its relative homology class. Due to the assumptions, this allows us to find a Seifert surface $\Sigma$ for $K$, homologous to $S$, with $\sel(K,\Sigma)=2g(\Sigma)-1=-\chi(\Sigma)$. So by Theorem~\ref{t:ribbon}, $K$ is the boundary of a Legendrian ribbon representing the class $[S]\in H_2(M,K)$.

Conversely, if $K=\p\Sigma_\Gamma$ for a Legendrian graph $\Gamma$ with $[\Sigma_\Gamma]=[S]\in H_2(M,K)$, then Lemma~\ref{l:euler} yields $\sel(K,S)=-\chi(\Sigma_\Gamma)=2 g(\Sigma_\Gamma)-1\geq 2g(S)-1$ because $S$ minimizes the genus in its relative homology class.

The last statement holds because any topological knot type in an overtwisted contact 3--manifold can be realized by a transverse knot $K$ in the complement of an overtwisted disk, and by repeated destabilization we can arrange $\sel(K,S)\geq 2g(S)-1$. 
\end{proof}

As it is shown in \cite{BaM}, the situation in tight contact structures is quite different: A transverse knot bounding a Legendrian ribbon is quasipositive in the sense of Rudolph~\cite{Ru}. Quasipositivity is quite a strong condition, as it implies chirality. For example, the figure--8 knot is not quasipositive since it is achiral, i.e. topologically equivalent to its mirror image. A classification of quasipositive knots up to 10 crossings is given in \cite{Ba}.

To conclude this note, we give a characterization of tightness in terms of Legendrian ribbons. 

\begin{thm}
A contact structure $\xi$ on a 3--manifold $M$ is tight if and only if every Legendrian ribbon has minimal genus among all embedded surfaces with the same boundary and in the same relative homology class. 
\end{thm}

\begin{proof}
Assume first that $\xi$ is tight and $\Sigma_\Gamma$ is a Legendrian ribbon with boundary $\p\Sigma_\Gamma=K$. By Lemma~\ref{l:euler} we have $\sel(K,\Sigma_\Gamma)=-\chi(\Sigma_\Gamma)$, and by Theorem~\ref{t:bennequin} every other Seifert surface $\Sigma$ for $K$ homologous to $\Sigma_\Gamma$ satisfies $\sel(K,\Sigma)\leq-\chi(\Sigma)$, so $g(\Sigma)\geq g(\Sigma_\Gamma)$.  

Conversely, if $\xi$ is overtwisted, then we construct a Legendrian ribbon $\Sigma_\Gamma$ of genus one whose boundary $K=\p\Sigma_\Gamma$ is the topological unknot, as shown in Figure~3.

\begin{figure}[ht]
\center
\includegraphics{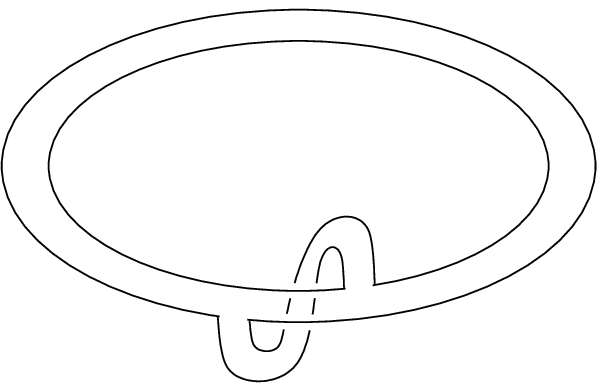}
\caption{}
\end{figure}

Since $\Sigma_\Gamma$ is contained in the neighbourhood of an overtwisted disc, i.e. in a ball, it is homologous rel $K$ to a disk and therefore not genus minimizing. 
\end{proof}


\begin{thebibliography}{12345}

\bibitem{Ba} S.~Baader, {\em Slice and Gordian numbers of track knots}, Osaka J. Math.  42 (2005), no. 1, 257--271.

\bibitem{BaM} S.~Baader, M.~Ishikawa, {\em Legendrian graphs and quasipositive digrams}, arXiv:math.GT/0609592. 

\bibitem{Be} D.~Bennequin, {\em Entrelacements et equations de Pfaff}, Ast{\'e}risque 107--108 (1983), 83--161. 

\bibitem{Dy} K.~Dymara, {\em Legendrian knots in overtwisted contact strucures}, arXiv:math.GT/0410122. 

\bibitem{El} Y.~Eliashberg, {\em Contact $3$--manifolds twenty years since J.~Martinet's work}, Ann. Inst. Fourier 42, 1--2 (1992), 165--192. 

\bibitem{Et} J.~Etnyre, {\em Legendrian and transversal knots}, preprint 2003. 

\bibitem{FT} D.~Fuchs, S.~Tabachnikov, {\em Invariants of Legendrian and transverse knots in the standard contact space}, Topology 36 no. 5 (1997), 1025--1053. 

\bibitem{Ka} Y.~Kanda, {\em On the Thurston--Bennequin invariant of Legendrian knots and nonexactness of Bennequin's inequality}, Invent. Math. 133 no. 2 (1998), 227--242. 

\bibitem{Ru} L.~Rudolph, {\em Constructions of quasipositive knots and links III. A characterization of quasipositive Seifert surfaces}, Topology 31 no. 2 (1992), 231--237. 

\end{thebibliography}
\end{document}